\documentclass[11pt, onecolumn]{IEEEtran}
\usepackage{cite}
\usepackage{amsmath,amssymb,amsfonts}
\usepackage{algorithmic}
\usepackage{graphicx}
\usepackage{textcomp}
\usepackage{euscript}
\usepackage{latexsym,,color,amstext,bm}
\usepackage{tikz}
\usepackage{circuitikz}

\def\BibTeX{{\rm B\kern-.05em{\sc i\kern-.025em b}\kern-.08em
    T\kern-.1667em\lower.7ex\hbox{E}\kern-.125emX}}
\usepackage{amsmath,amssymb,euscript,yfonts,psfrag,latexsym,dsfont,graphicx}
\usepackage{bbm,color,amstext,wasysym,subfig,flushend,parskip,balance}
\usepackage{algorithm}
\usepackage{algorithmic}
\graphicspath{{./},{./figures/}}

\usepackage{tikz}
\usetikzlibrary{fit,positioning}

\newtheorem{theorem}{Theorem}

\newtheorem{prop}{Proposition}
\newtheorem{problem}[theorem]{Problem}

\newtheorem{definition}[theorem]{Definition}

\newcommand{\mR}{{\mathbb R}}

\newcommand{\cP}{{\mathcal P}}

\newcommand{\cV}{{\mathcal V}}

\newcommand{\one}{{\mathds 1}}

\newcommand{\bC}{{\mathbf C}}
\newcommand{\bM}{{\mathbf M}}
\newcommand{\bQ}{{\mathbf Q}}
\newcommand{\support}{{\rm Supp}}

\newcommand{\diag}{\operatorname{diag}}

\newcommand{\mfp}{\mathfrak{p}}

\definecolor{grey}{rgb}{0.6,0.6,0.6}
\definecolor{lightgray}{rgb}{0.97,.99,0.99}

\usepackage{soul}
\begin{document}
\title{
Control and estimation of multi-commodity network flow under aggregation
}

\author{Yongxin Chen, Tryphon T.\ Georgiou and Michele Pavon
\thanks{Y.\ Chen is with the School of Aerospace Engineering, Georgia Institute of Technology, Atlanta, GA 30332, USA;yongchen@gatech.edu }
\thanks{T.T.\ Georgiou is with the Department of Mechanical and Aerospace Engineering, University of California, Irvine, CA 92697, USA; tryphon@uci.edu}
\thanks{M. Pavon is with the Division of Science, New York University Abu Dhabi, U.A.E.; michele.pavon@nyu.edu}
\thanks{This work was supported in part by the NSF under grants 1942523 and 2206576, the AFOSR under FA9550-23-1-0096,  the ARO under W911NF-22-1-0292, and the NYUAD under grant 76/71260/ADHPG.}}
\maketitle

\begin{abstract}

A paradigm put forth by E. Schr\"odinger in 1931/32, known as Schr\"odinger bridges, represents a formalism to pose and solve control and estimation problems seeking a perturbation from an initial control schedule (in the case of control), or from a prior probability law (in the case of estimation), sufficient to reconcile data in the form of marginal distributions and minimal in the sense of relative entropy to the prior. In the same spirit, we consider traffic-flow and apply a Schr\"odinger-type dictum, to perturb minimally with respect to a suitable relative entropy functional a prior schedule/law so as to reconcile the traffic flow with scarce aggregate distributions on families of indistinguishable individuals. Specifically, we consider the problem to regulate/estimate multi-commodity network flow rates based only on empirical distributions of commodities being transported (e.g., types of vehicles through a network, in motion) at two given times. Thus, building on Schr\"odinger's large deviation rationale, we develop a method to identify {\em the most likely flow rates (traffic flow)}, given prior information and aggregate observations. 
Our method further extends the Schr\"odinger bridge formalism to the multi-commodity setting, allowing commodities to exit or enter the flow field as well (e.g., vehicles to enter and stop and park) at any time. The behavior of entering or exiting the flow field, by commodities or vehicles, is modeled by a Markov chains with killing and creation states. Our method is illustrated with a numerical experiment.
\end{abstract}

\section{Introduction}
{\em Inverse Problems} constitute a large class of typically ill-posed problems of central importance in all branches of science. In an inverse problem, one seeks to derive a model (a function, a field, a probability distribution, etc.) from a set of observations. Examples are ubiquitous in system identification, spectral estimation, computed tomography, deconvolution, inverse scattering, weather prediction, and so on. For instance,  image  deblurring may be viewed as a deconvolution problem in the plane, where resolution is to be restored based on priors on features of objects and texture.

{\em Regularization} is a process in mathematics, statistics, and machine learning, that consists in adding information to solve ill-posed problems and/or to prevent {\em overfitting}\footnote{An overfitted model is a statistical model that contains more parameters than can be justified by the data thereby violating {\em Novacula Occami}: ``Frustra fit per plura quod potest fieri pauciora". Novacula Occami should perhaps be better translated ``Ockham's comb" rather than Ockham's razor as it is customary.}. For instance, in compressed sensing, sparsity of the solution is the added information. Other times, it consists in imposing smoothness of the solution or penalazing some norm of the solution. From a Bayesian point of view, many regularization techniques correspond to imposing certain prior distributions on model parameters.  Tikhonov's regularization, for instance, is widely  used in finite and infinite-dimensional context,
see e.g., \cite{cucker2002mathematical}. 

A most powerful  paradigm ({\em entropy regularization}) to learn a probability distribution from scarce information was put forward by Ludwig Boltzmann \cite{Boltzmann1877beziehung} in 1877 and by Erwin Schr\"{o}dinger \cite{Sch31,Sch32} in 1931/32. The {\em most likely probability distribution} can be characterized as the solution of a {\em maximum entropy problem.} This inference method has, in the meantime, proved very fruitful in several branches of science: We mention the far reaching work of Jaynes, Burg, Dempster and Csisz\'{a}r \cite{jaynes1982rationale,burg1975maximum,dempster1972covariance,csiszar1975divergence}. In the case of a discrete version of the dynamical problem considered by Schr\"{o}dinger, this method can be fruitfully reformulated as a {\em Markov decision problem}. We show in this paper that a suitable modification of this inference method can be applied to characterize  network flows under very meager information. We consider indeed multi-commodity  flows  where only aggregate information is available at some initial and final times. Moreover, some of the agents may enter or exit the flow during the time interval of interest. Hence, the total number of travelling vehicles, packets, etc. is in general not preserved.

To give a more precise formulation of the problem, consider for instance a network of roads and highways. At some initial time $t=0$, we observe the distribution of different types of vehicles traveling on the network. The distribution $\mu_0$ of the various classes (cars, trucks, etc.) is  defined on the vertices of the network (e.g., at crossroads, intersections, or cities, for such a network). At some final time $t=N$, we observe a similar distribution $\mu_N$. We seek to determine the most likely network flow on the discrete time interval $[0,N]$ which is consistent with the available scant information. An effective framework to study such problems is a suitable extension   of what is called {\em regularized optimal transport}, also known as {\em discrete Schr\"odinger Bridges} \cite{PavTic10,Cut13,BenCarCut15,GeoPav14,peyre2019computational,haasler2020optimal,chen2021stochastic,chen2021controlling,chen2021optimal2,koehl2021physics,zhou2021optimal,arque2022approximate,leleux2022design,di2022maximal,courtain2023relative}. We show that such problems admit a formulation as suitable Markov decision problems in the spirit of   \cite[Section 6]{chen2021optimal2}.
We are actually able to solve the problem in the more complex situation, as noted, when some of the vehicles enter or exit the flow during $[0,N]$.

In a similar spirit, in the multi-commodity setting that we study herein, and in the context of transportation, we assume knowledge of the fraction of vehicles that are private cars, taxis, trucks, and so on (that more generally may be thought of as different commodities, species, etc.). We also assume knowledge of a {\em prior} probability on the flow across the network for each specific group of vehicles; such information can in principle be provided by historic data. With such data and assumptions in place, we are interested in identifying the most likely path that each commodity has followed while being transported across the network. A similar problem may be formulated when commodities have sources and sinks. 

In the present work we show that by suitably extending the theory of Schr\"odinger Bridges it is possible to answer such questions in spite of the scarce, aggregate information available at the initial and final times. As a second contribution of the paper, we address the case where vehicles transporting commodities may enter or exit the flow during the prescribed window of time. The sudden appearance or disappearance of vehicles can be modeled probabilistically, via notions referred to as creation or killing, and can be derived along the lines of  \cite{chen2021most} that dealt with single-commodity transport with killing, in continuous-time.

The paper is structured as follows. In Section \ref{GBProblems}, we briefly recall the formulation and key results for the single commodity mass-preserving Schr\"odinger Bridge problem. Section \ref{MCNF} formulates and solves the multi-commodity network flow problem. In Section \ref{CRKIL}, we study the same problem in the presence of creation and killing.



\section{Single-commodity traffic flow over networks}\label{GBProblems}

We begin by discussing a paradigm of great significance in single-commodity network flows. It amounts to designing probabilistic transitions between nodes, and thereby probability laws on path spaces, so as to reconcile marginal distributions with priors that reflect the structure of the network and objectives on transference of resources across the network. More precisely, we wish to study a generalization of the so called discrete-time Schr\"odinger bridge problem. This dynamic formulation echoes the fluid dynamic problem associated to the classical (continuous time and space) Schr\"odinger bridge problem.

Consider a directed, strongly connected, aperiodic graph ${\bf G}=(\cV,\mathcal E)$ with vertex set $\cV=\{1,2,\ldots,n\}$ and edge set $\mathcal E\subseteq \cV\times\cV$. Consider trajectories/paths on this graph over the time set ${\mathcal T}=\{0,1,\ldots,N\}$. 
We seek a probability distribution on the space of paths $\cV^{N+1}$ with prescribed initial and final marginals
$\mu_0$ and $\mu_N$, respectively, and such that the resulting random evolution
is closest to a ``prior'' measure $\bQ$ in a suitable sense.

The prior law for our problem is a Markovian evolution with transition kernel $A$, which may be assumed to be time-homogenous for simplicity. 
In accordance with 
the topology of the graph, $A(i,j)=0$ whenever $(i,j)\not\in\mathcal E$.  
We assume that $\mu_0$ is positive  on $\cV$, i.e.
\begin{equation}\label{eq:mupositive}
\mu_0(i)>0\mbox{ for all }i\in\cV.
\end{equation}
The Markovian kernel $A$, together with the measure $\mu_0(\cdot)$, 
 induces
a measure $\bQ$ on the space of trajectories, which assigns to a path  $(i_0,i_1,\ldots,i_N)\in \cV^{N+1}$ the value
\begin{equation}\label{prior}
\bQ(i_0,i_1,\ldots,i_N)=\nu_0(i_0)A(i_0,i_1)\cdots A(i_{N-1},i_N).
\end{equation}
\begin{definition} We denote by ${\mathcal P}(\mu_0,\mu_N)$ the family of probability distributions on $\cV^{N+1}$ having the prescribed marginals $\mu_0$ and $\mu_N$.
\end{definition}

We seek a distribution in this set which is closest to the prior $\bQ$ in {\em relative entropy} (divergence, Kullback-Leibler index)  defined by 
\begin{equation*}
{\rm KL}(\bM \|\bQ):=\left\{\begin{array}{ll} \sum_{x}\mathfrak \bM (x)\log\frac{\bM (x)}{\bQ (x)}, & \support (\bM)\subseteq \support (\bQ),\\
+\infty , & \support (\bM)\not\subseteq \support (\bQ).\end{array}\right.
\end{equation*}
Here, by definition,  $0\cdot\log 0=0$.
This brings us to the 
so-called
{\em Schr\"odinger Bridge  Problem} (SBP):

\begin{problem}\label{prob:optimization}
Determine
 \begin{eqnarray}\label{eq:optimization}
\bM^\star={\rm argmin}\{ {\rm KL}(\bM\|\bQ) \mid  \bM\in {\mathcal P}(\mu_0,\mu_N)
\}.
\end{eqnarray}
\end{problem}

We parameterize $\bM$ as
\begin{equation}\label{markovianmeasure}
\bM(i_0,i_1,\ldots,i_N)=\mu_0(i_0)\pi_{i_0i_1}(0)\cdots \pi_{i_{N-1}i_N},
\end{equation}
and let $p_t$ be the one-time marginals of $\bM$, i.e.,
\[p_t(i_t) = \sum_{i_{\ell\neq t}}\bM(i_0,i_1,\ldots,i_N), \quad t\in\mathcal T.\] 
We finally have the update mechanism
\begin{equation}\label{update}
p_{t+1}(i_{t+1})=\sum_{i_t\in\cV} p_t(i_t) \pi_{i_{t}i_{t+1}}(t)
\end{equation}
which, in (column) vector form, is 
\begin{equation}
\label{eq:fp}
p_{t+1}=\Pi'(t) p_{t}.
\end{equation}
Here $\Pi=\left[ \pi_{ij}(t)\right]_{i,j=1}^n$ is the transition matrix and prime denotes transposition.
Using (\ref{prior})-(\ref{markovianmeasure}) we obtain
\[{\rm KL}(\bM\| \bQ)={\rm KL}(\mu_0\|\nu_0) + \sum_{t=0}^{N-1}\sum_{i_t}{\rm KL}
(\pi(t)\|A(t))p_t(i_t).
\]
 We have the following theorem \cite{PavTic10,GeoPav14,chen2016robust}:

\begin{theorem} \label{maincl}Suppose there exists a pair of nonnegative functions $(\varphi,\hat{\varphi})$ defined  on $\{0,1,\ldots,N\}\times{\cal V}$ and satisfying the system
\begin{eqnarray}\label{scsi1}\varphi(t,i_t)=\sum_{i_{t+1}}A_{i_{t}i_{t+1}}(t)\varphi(t+1,i_{t+1}),\\
\hat{\varphi}(t+1,i_{t+1})=\sum_{i_t}  A_{i_{t}i_{t+1}}(t)\hat{\varphi}(t,i),\label{scsi2}
\end{eqnarray}
for $t=0,1,\ldots, N-1$, as well as the boundary conditions
\begin{equation}\label{scsi3}
\varphi(0,i_0)\hat{\varphi}(0,i_0)=\mu_0(i_0),\; \varphi(N,i_N)\hat{\varphi}(N,i_N)=\mu_N(i_N),
\end{equation}
for $i_t\in\cV$ and $t\in\{0,N\}$, accordingly.
Suppose moreover that $\varphi(t,i)>0,\; \forall t=0, 1, \ldots, N,\, \forall i\in{\cal V}$. Then, the Markov distribution $\bM^\star$ in ${\mathcal P}(\mu_0,\mu_N)$ having transition probabilities
\begin{equation}\pi^\star_{i_{t}i_{t+1}}(t)=A_{i_{t}i_{t+1}}(t)\frac{\varphi(t+1,i_{t+1})}{\varphi(t,i_{t})}
\end{equation}
solves Problem \ref{prob:optimization}. 
\end{theorem}
Notice that if $(\varphi,\hat{\varphi})$ satisfy (\ref{scsi1})-(\ref{scsi2})-(\ref{scsi3}), so does the pair $(c\varphi,\frac{1}{c}\hat{\varphi})$ for all $c>0$. Hence, uniqueness for the Schr\"{o}dinger system is always intended as uniqueness of the ray.
Under the assumption that  the entries of the matrix product $A^N$
are all positive, there exists a (unique in the sense of ray) solution to the system (\ref{scsi1})-(\ref{scsi3}) which can be computed through a Fortet-IPF-Sinkhorn iteration \cite{For40,deming1940least,sinkhorn1964relationship,Cut13,GeoPav14}. The infinite-horizon counterpart of it has been considered in \cite{SaeYen09,chen2021optimal,chen2021fast}.

\section{Multi-commodity traffic flow}\label{MCNF}
Consider a traffic flow consisting of $K$ commodities. Each commodity may correspond to one origin-destination pair. Other ways to differentiate commodities may include vehicle types. Now assume that the prior population distribution of the $K$ commodities follows the probability (vector) $\mfp = [\mfp_1,\mfp_2,\ldots,\mfp_K]^T$, that is,  $\mfp_k$ represents the portion of the population corresponding to the $k$-th commodity. Suppose each individual in the $k$-th commodity follows a time-homogenous Markov model with transition kernel $A_{k}$ and initial distribution $r_k$. The transition kernel for each commodity can be time-varying, however, we henceforth use time-homogenous kernel to keep the notation simple. 

The above defines a probabilistic model for the traffic. It is a  hierarchical  model, where the first level is the commodity type which follows the distribution $\mfp$. The second level represents the dynamics over the traffic network captured by the transition kernel $A_{k}$ with initial distribution $r_k$, for each commodity. A traffic flow then can be viewed as a collection of independent samples from this probabilistic model parametrized by the tuple $(\mfp, A_1, r_1,\ldots, A_K, r_K)$.

We are interested in the inference problem of estimating the group transport of the commodities with limited and aggregated data. In the aggregate measurements, the individuals from different commodities are indistinguishable. Thus, the measurement at any time-instance is a histogram, representing the distribution of the total population over the traffic network. The measurements are limited in the sense that we can only make measurement every once in a while. In particular, we consider the setting where the traffic flow is only measured at two points in time, $t=0$ and $t=N$. 

Suppose the population measurement is $\mu_0$ at the initial time $t=0$ and $\mu_N$ at terminal time $t=N$. Here the measurements 
$\mu_0, \mu_N$ are normalized to be probability vectors dividing by the total number $L$ of vehicles. Our goal is to infer the most likely evolution of each commodity. More precisely, we want to recover the most likely distribution on the space of trajectories over the traffic network that the vehicles may have takes so as to match the measurements $\mu_0, \mu_N$, taken at $t=0, N$. When there is only a single commodity, i.e., $K=1$, the problem clearly reduces to a standard Schr\"odinger bridge problem. Thus, this problem represents a generalization of the standard SBP to the multi-commodity setting. 

To describe the evolution of the group behavior of the $K$ commodities, let $\bM^k\in\mR^{n^{(N+1)}}$ be the normalized population distribution of the $k$-th commodity over the path space $\cV^{N+1}$. It means that the portion of the population in the $k$-th commodity travels along the graph path $i_0, i_1, \ldots, i_N$ is then $\bM^k_{i_0,i_1,\ldots,i_N}$. Note that due to normalization, $\sum_{k=1}^K |\bM^k| = 1$ where $|\bM^k|$ denotes the 1-norm of the tensor $\bM^k$, describing the portion of population that belongs to the $k$-th commodity. The tensors $\bM^1, \bM^2, \ldots, \bM^K$ fully characterize the group behavior of the whole population. Our goal is to find the most likely tuple $\bM^1, \bM^2, \ldots, \bM^K$ given the observations $\mu_0, \mu_N$. It turns out that this multi-commodity network flow also satisfies the large deviation principle as in the following theorem.
\begin{theorem}
The probability of an empirical population distribution $(\bM^1,\bM^2, \cdots,\bM^K)$ in the traffic flow, parameterized by the tuple $(\mfp, A_1, r_1,\ldots, A_K, r_K)$, is
	\begin{equation}
		{\rm Prob} (\bM^1,\bM^2, \cdots,\bM^K) \approx \exp[-L \times I(\bM^1,\bM^2, \cdots,\bM^K)],
	\end{equation}
where $L$ is the total number of the population, and
\begin{align}\label{eq:rate}
		I(\bM^1,\bM^2, \cdots,\bM^K) &=  \sum_{k=1}^K {\rm KL} (\bM^k \|\,| 
\bM^k|r_k(i_0)A_k(i_0,i_1)\cdots \\
  & \hspace*{-15pt}\cdots A_k(i_{N-1},i_N)) + \sum_{k=1}^K |\bM^k|\log \frac{|\bM^k|}{p_k}.\nonumber
	\end{align}
 is the rate function.
\end{theorem}
{\bf Proof}
The empirical distribution $(\bM^1,\bM^2, \cdots,\bM^K)$ assigns probability $\bM^k_{i_0,i_1,\ldots,i_N}$ to the observation $k, i_0,i_1,\ldots,i_N$. The prior model assigns probability $\mfp_k r_k (i_0) A_k(i_0,i_1)\cdots A_k(i_{N-1}, i_N)$ to the same observation. By the standard theory of large deviation, the rate function $I$ is the KL divergence between these two probability distributions. The expression \eqref{eq:rate} then follows a straightforward calculation of the KL divergence. 

\begin{figure}[tb]
	\centering
 \includegraphics[width=0.41\textwidth]{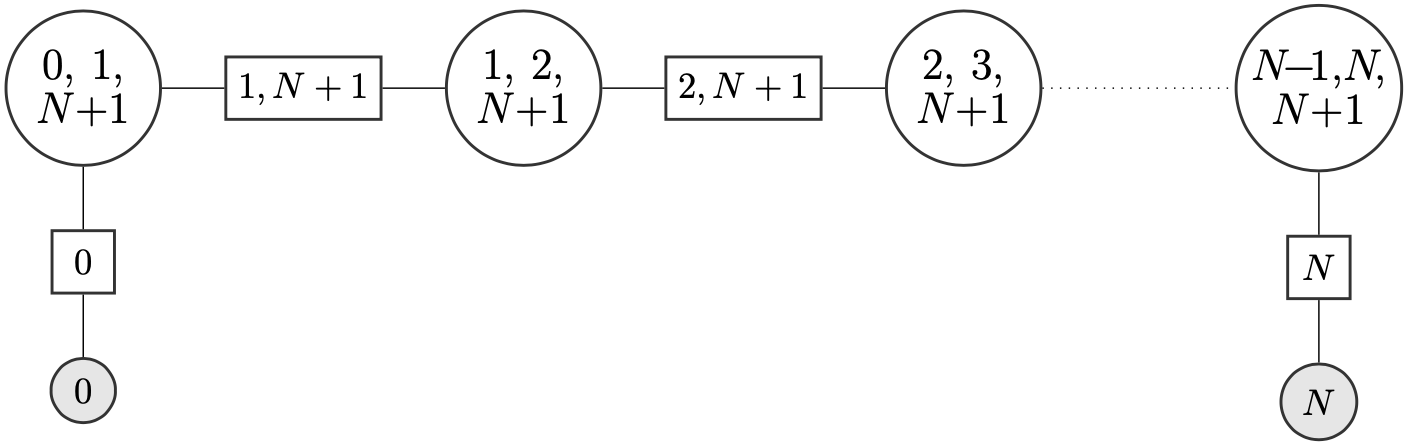}
	\caption{Junction tree for the graphical OT \eqref{eq:graphOT}}
	\label{fig:traffic_graph}
\end{figure}

In light of the above large deviation result, our problem of recovering the most likely evolution can be formulated as the optimization
	\begin{subequations}\label{eq:MCTF}
	\begin{align}\nonumber
\min_{\bM^1,\cdots,\bM^K} & \sum_{k=1}^K  {\rm KL} (\bM^k \|\,| \bM^k|r_k(i_0)A_k(i_0,i_1)\cdots\\
& A_k(i_{N-1},i_N)) + \sum_k |\bM^k|\log \frac{|\bM^k|}{p_k}
		\\\label{eq:13b}
 \mbox{ subject to }\hspace*{10pt} & \sum_k \sum_{i_1,i_2,\cdots,i_N} \bM^k_{i_0,i_1,\ldots,i_N} = \mu_0(i_0),\\\label{eq:13c}
  & \sum_k \sum_{i_0,i_1,\cdots,i_{N-1}} \bM^k_{i_0,i_1,\ldots,i_N} = \mu_N(i_N).
	\end{align}
	\end{subequations}
Just like the standard SBP, our inference problem can be formulated as an optimal transport problem, albeit multi-marginal, with entropic regularization. To this end, let $\bM = [\bM^1,\bM^2,\cdots,\bM^K]$. Then, straightforward calculations yield the following reformulation of \eqref{eq:MCTF}. 
\begin{prop}
The multi-commodity traffic flow inference problem \eqref{eq:MCTF} can be reformulated as a multi-marginal optimal transport with entropy regularization
	\begin{subequations}\label{eq:graphOT}
\begin{align}\label{eq:graphOT1}
		\min_{\bM} & \langle \bM, \bC\rangle+ \langle \bM, \log \bM\rangle 
		\\\nonumber
&\mbox{ subject to }\hspace*{5pt}  (\ref{eq:13b}-\ref{eq:13c})
	\end{align}
where the transport cost tensor is 
\begin{align}\nonumber
\bC_{k,i_0,i_1,\cdots,i_N} &= -\log r_k(i_0)A_k(i_0,i_1)\cdots\\&\cdots\label{eq:graphOT3} A_k(i_{N-1},i_N) -\log p_k.
	\end{align}
	\end{subequations}
\end{prop}
The 	above entropy-regularized, multi-marginal optimal transport (MOT) problem can be solved using the standard Sinkhorn algorithm. The Sinkhorn algorithm for \eqref{eq:graphOT} is a block ascent algorithm for its dual \cite{HaaSinChe20}. Each iteration requires a projection operation to compute a marginal distribution of the tensor $\bM$. Thus, a generic Sinkhorn solver of this type has computational complexity that scales exponentially with the number of marginals (in our case, $N+2$). Fortunately, the optimization \eqref{eq:graphOT} is a MOT problem with graph-structured cost \cite{HaaSinChe20}. In particular, the cost tensor in \eqref{eq:graphOT3} can be written as
	\begin{equation}\label{eq:Ctensor}
		\bC_{k,i_0,i_1,\cdots,i_N} = \sum_{t=1}^N \bC^{t}_{k,i_{t-1},i_t} 
	\end{equation}
where
	\begin{equation*}
		\bC^1_{k,i_0,i_1} = -\log r_k(i_0)A_k(i_0,i_1) - \log p_k, 
	\end{equation*}
and
	\begin{equation*}
		\bC^t_{k,i_{t-1},i_t} = -\log A_k(i_{t-1},i_t), \quad 2\le t\le N.
	\end{equation*}
We remark that,  even though the total cost tensor $\bC$ is $N+2$ dimensional, it is the summation of $N$ tensors of dimension $3$. This decomposition enables us to exploit the graphical structure of the cost tensor to greatly reduce the computational complexity of the Sinkhorn algorithm. The cost in \eqref{eq:Ctensor} is associated with the junction tree in Figure \ref{fig:traffic_graph}.
This junction tree is associated with a graph with $N+2$ nodes. The nodes $0, 1, \ldots, N$ correspond to the vehicle distributions of the total population at each time point. The node $N+2$ is used to model the mass distribution over different commodities. 

\begin{theorem}
The solution to \eqref{eq:graphOT} is characterized by the system of equations
	\begin{subequations}\label{eq:FB}
	\begin{align}
		\varphi(t,i_t, k) &= \sum_{i_{t+1}} A_k(i_t,i_{t+1}) \varphi(t+1,i_{t+1}, k), 
		\label{eq:FB1}
		\\
		\hat\varphi(t+1,i_{t+1},k) &= \sum_{i_{t}} A_k(i_t,i_{t+1}) \hat\varphi(t,i_{t}, k),
		\label{eq:FB2}
		\\
		\varphi(0,i_0) &= \sum_{i_1,k}r_k(i_0)p_kA_k(i_0,i_{1}) \varphi(1,i_{1}, k), \label{eq:FB3}
		\\
		\hat\varphi(1,i_1,k) &= \sum_{i_0} r_k(i_0)p_kA_k(i_0,i_1) \hat\varphi(0,i_0), \label{eq:FB4}
		\\
		\varphi(N-1,i_{N-1},k) &= \sum_{i_{N}} A_k(i_{N-1},i_{N}) \varphi(N,i_{N}), \label{eq:FB5}
		\\
		\hat\varphi(N,i_N) &= \sum_{i_{N-1},k}A_k(i_{N-1},i_{N}) \hat\varphi(N-1,i_{N-1}, k), \label{eq:FB6}
	\end{align}
	\end{subequations}
for $t = 1, \ldots, N-2 $, with boundary conditions
	\begin{equation}\label{eq:BD}
		\varphi(0, \cdot) \hat\varphi(0,\cdot) = \mu_0,\quad \varphi(N,\cdot)\hat\varphi(N,\cdot) = \mu_N.
	\end{equation}
Moreover, the transition probabilities of the solution are
	\begin{eqnarray*}
		\pi_t^k (i_t,i_{t+1}) &=& A_k(i_t,i_{t+1}) \frac{\varphi(t+1,i_{t+1},k)}{\varphi(t,i_t,k)},
		\\
		\pi_0^k(i_0,i_1) &=& A_k(i_0,i_{1}) \frac{\varphi(1,i_{1},k)}{\sum_{i_1}A_k(i_0,i_{1})\varphi(1,i_{1},k)}
		\\
		\pi_{N-1}^k(i_{N-1},i_N) &=& A_k(i_{N-1},i_{N}) \frac{\varphi(N,i_{N})}{\varphi(N-1,i_{N-1},k)}
	\end{eqnarray*}
and the marginal distributions are 
	\begin{eqnarray*}
		\mu_t^k(i_t) &=& \varphi(t,i_t, k) \hat \varphi(t,i_t,k), \quad t= 1, \ldots, N-1
		\\
		\mu_0^k(i_0)&=& \hat\varphi(0,i_0) \sum_{i_1}r_k(i_0)p_kA_k(i_0,i_{1}) \varphi(1,i_{1}, k)
		\\
		\mu_N^k(i_N)&=& \varphi(N,i_N) \sum_{i_{N-1}} A_k(i_{N-1},i_{N}) \hat\varphi(N-1,i_{N-1}, k).
	\end{eqnarray*}
\end{theorem}

The graphical OT can be efficiently solved with the Sinkhorn belief propagation algorithm \cite{SinHaaChe2022,HaaSinChe20,HaaRinCheKar21}, which is a combination of the Sinkhorn algorithm for OT and the belief propagation algorithm for probabilistic graphical models. Unlike the vanilla Sinkhorn  \cite{sinkhorn1964relationship}, whose complexity scales exponentially with the number of marginals, the complexity of the Sinkhorn belief propagation algorithm is determined by the largest node degree  in the junction tree associated with the cost tensor. For our multi-species problem \eqref{eq:graphOT}, the Sinkhorn belief propagation algorithm is specialized to Algorithm \ref{alg:forward_backward}.

\begin{algorithm}[ht]
   \caption{Sinkhorn Belief Propagation Algorithm for Multi-commodity Traffic Flow}
   \label{alg:forward_backward}
\begin{algorithmic}
   \STATE Initialization: $\hat\varphi(0,\cdot) = {\bf 1}$ 
   \WHILE{not converged}
   \STATE \textbf{Forward pass:}
   \FOR{$t = 1,2,\ldots,N-1$}
        \STATE Update  $\hat\varphi(t, \cdot,\cdot)$ using \eqref{eq:FB2} or \eqref{eq:FB4}
    \ENDFOR
    	 \STATE Update  $\hat\varphi(N, \cdot)$ using \eqref{eq:FB6}
	 \STATE Update  $\varphi(N, \cdot)$ using \eqref{eq:BD}
    \STATE \textbf{Backward pass:}	
    \FOR{$t = N-1,\ldots,1$}
        \STATE i) Update  $\varphi(t,\cdot,\cdot)$ using \eqref{eq:FB1} or \eqref{eq:FB5}
    \ENDFOR
    	\STATE Update  $\varphi(0, \cdot)$ using \eqref{eq:FB3}
	\STATE Update  $\hat\varphi(0, \cdot)$ using \eqref{eq:BD}
    \ENDWHILE
\end{algorithmic}
\end{algorithm}


\section{Most likely network flow with creation and killing}\label{CRKIL}
We now turn our attention to the traffic-flow estimation in the case where commodities may disappear or be suddenly introduced, at various times and locations, during transport.  The appearance and disappearance are referred to as creation and killing. 
In applications as in traffic flow, the creation may model the situation where a vehicle enters the traffic flow from a parking state, and the killing models the converse.  We focus on the single-commodity case, but the same idea can be applied to the multi-commodity setting in a completely analogous manner.

In the single-commodity case, the prior dynamics is encoded in a transition kernel $A$ over the traffic network ${\bf G}=(\cV,\mathcal E)$ as before, but which however, may not be necessarily a stochastic matrix in the present context, as commodities are created or destroyed. That is, whether time-homogeneous or not,
	\begin{equation}
		A \mathds{1} \neq \mathds{1}
	\end{equation}
in general. In fact, we assume that the row sums are less than $1$. This implies that, on each node $i \in \cV$ in the graph, the random walk has probability of $1-\sum_{j} A (i,j)$ to be killed/absorbed. In our traffic-flow estimation problem, we observe two marginal distribution $\mu_0$ at time $0$ and $\mu_N$ at time $N$ that may not have the same mass, in general. Our goal is once again to estimate the most likely traffic flow, given the prior dynamics, that is consistent with the two observed marginals. 

To account for the unbalance in mass between the marginals, we introduce a parking state for killing and creation. A similar idea has been used in \cite{chen2021most} to study an unbalanced Schr\"odinger bridge problem over a continuous state-space. We assume that the prior probability of going from the parking state to the graph nodes (creation rate at each node) is $c \in \mR_+^n$ with $c^T \one\le 1$. Augmenting the state space by the parking state, we obtain a Markov chain with transition kernel
	\begin{equation}
		\hat A = \left[\begin{matrix}A & b\\ c^T & d\end{matrix}\right]
	\end{equation}
where $b = \one-A \one$ and $d = 1- c^T \one$. We assume that the total number of vehicles in the augmented traffic network (including the parking state) is fixed. Without loss of generality, suppose that the marginals $\mu_0$ and $\mu_N$ have been normalized with respect to this number, that is, $\mu_0^T \one \le 1, \mu_N^T \one \le 1$. We define the augmented marginal distributions 
	\begin{subequations}
	\begin{equation}
		\hat\mu_0 = \left[\begin{matrix}\mu_0\\ 1 - \mu_0^T \one\end{matrix}\right]
	\end{equation}
and
	\begin{equation}
		\hat\mu_N = \left[\begin{matrix}\mu_N\\ 1 - \mu_N^T \one\end{matrix}\right].
	\end{equation}
	\end{subequations}
Let $\bQ$ be the distribution over the trajectories induced by the prior dynamics. It assigns to a path  $(i_0,i_1,\ldots,i_N)\in \cV^{N+1}$ the value
\begin{equation}\label{priornew}
\bQ(i_0,i_1,\ldots,i_N)=\hat \nu_0(i_0)\hat A(i_0,i_1)\cdots \hat A(i_{N-1},i_N).
\end{equation}
Denote by $\cP(\hat \mu_0, \hat \mu_N)$ all the Markov chains over the augmented state space with marginal distribution $\hat \mu_0$ at time $0$ and marginal distribution $\hat \mu_N$ at time $N$. Then our estimation problem becomes a Schr\"odinger bridge problem as follows.

\begin{problem}\label{prob:secondmaxent}
 \begin{equation}\label{eq:secondmaxent}
\min\{ {\rm KL} (\bM \| \bQ) \,|\,  \bM\in {\mathcal P}(\hat\mu_0,\hat\mu_N)
\}
\end{equation}
\end{problem}

Applying the standard Schr\"odinger bridge theory (as in Theorem \ref{maincl}) we obtain the following characterization of the solution to \eqref{eq:secondmaxent}.
\begin{theorem}
The Schr\"odinger system
	\begin{subequations}\label{eq:SBsysCK}
	\begin{eqnarray}
		\left[\begin{matrix}
		A^T & c\\b^T & d
		\end{matrix}\right]
		\left[\begin{matrix}
		\hat\varphi_t \\ \hat\psi_t
		\end{matrix}\right]
		&=&
		\left[\begin{matrix}
		\hat\varphi_{t+1} \\ \hat\psi_{t+1}
		\end{matrix}\right]
		,\quad t = 0, 1,\ldots, N-1
		\\
		\left[\begin{matrix}
		A & b \\c^T & d
		\end{matrix}\right]
		\left[\begin{matrix}
		\varphi_{t+1} \\ \psi_{t+1}
		\end{matrix}\right]
		&=&
		\left[\begin{matrix}
		\varphi_{t} \\ \psi_{t}
		\end{matrix}\right]
		,\quad t = 0, 1,\ldots, N-1
		\\
		\left[\begin{matrix}
		\varphi_{0}\hat\varphi_0 \\ \psi_{0}\hat\psi_0
		\end{matrix}\right]
		&=&
		\left[\begin{matrix}
		\mu_0 \\ 1-\mu_0^T \mathds{1}
		\end{matrix}\right],
		\\
		\left[\begin{matrix}
		\varphi_{N}\hat\varphi_N \\ \psi_{N}\hat\psi_N
		\end{matrix}\right]
		&=&
		\left[\begin{matrix}
		\mu_N \\ 1-\mu_N^T \mathds{1}
		\end{matrix}\right]
	\end{eqnarray}
	\end{subequations}
admits a unique (up to a constant factor) solution. Moreover, the solution to the SBP problem \eqref{eq:secondmaxent} has transition matrix
	\begin{equation}\label{eq:D}
		\left[\begin{matrix}
		\diag(\varphi_t)^{-1} A \diag(\varphi_{t+1}) & \diag(\varphi_t)^{-1} b\psi_{t+1}
		\\ \frac{1}{\psi_t}c^T\diag(\varphi_t) & \frac{d\psi_{t+1}}{\psi_t}
		\end{matrix}\right]
	\end{equation}
with associated marginal flow
	\begin{equation}\label{eq:mu}
		\hat \mu_t = \left[\begin{matrix}
		\varphi_{t}\hat\varphi_t \\ \psi_{t}\hat\psi_t
		\end{matrix}\right].
	\end{equation}
\end{theorem}

Based on this result, we can recover the most likely evolution for the original Markov chain without the parking state. Its transition matrix is
	\begin{equation}
		\diag(\varphi_t)^{-1} A \diag(\varphi_{t+1}).
	\end{equation}
Note that, in general, this is not a stochastic matrix. This result should be compared with Theorem \ref{maincl} which does not involve killing or creation.


\section{Numerical Examples}
\begin{figure}[h]
	\centering
	\includegraphics[width=0.41\textwidth]{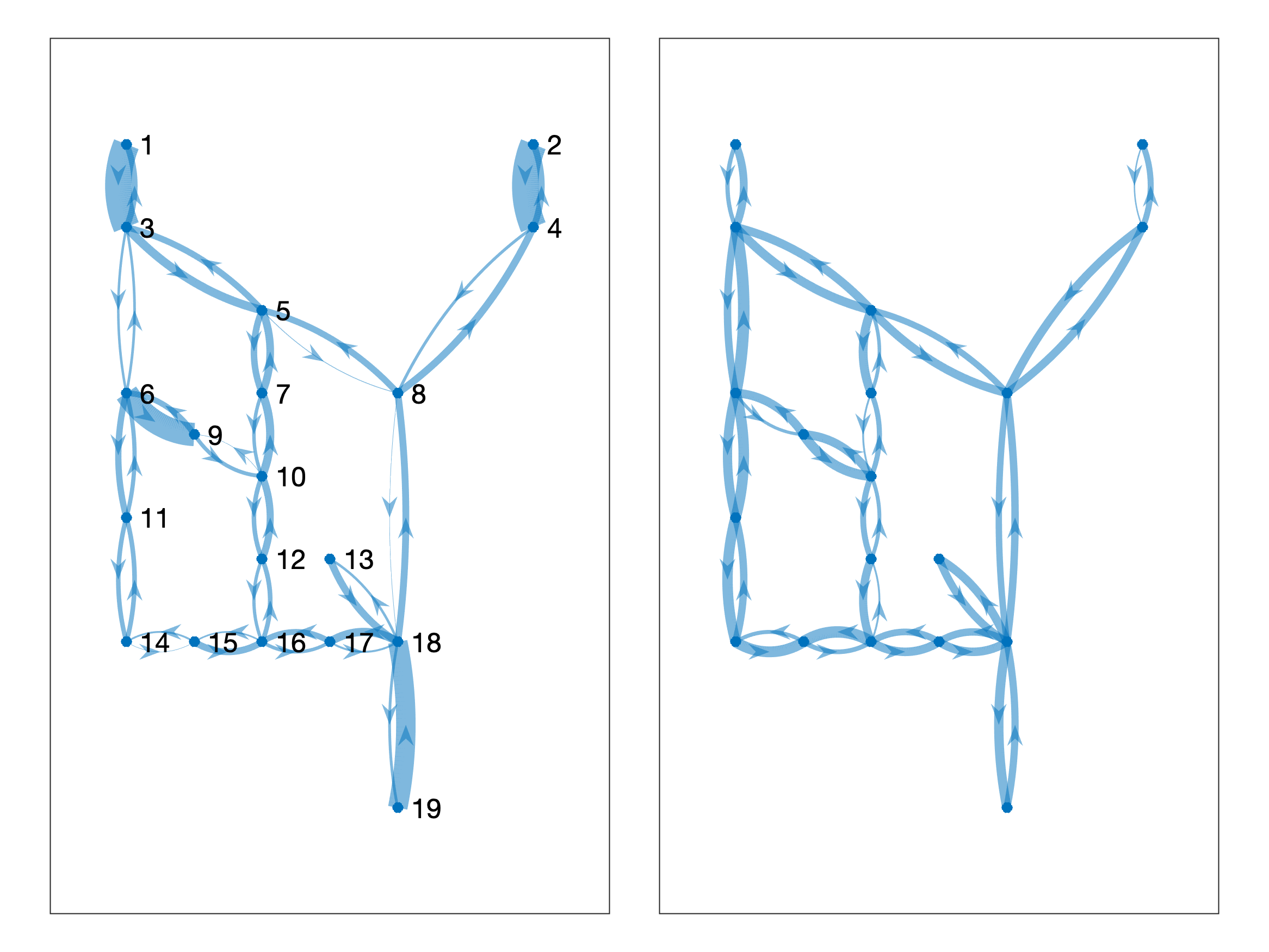}
	\caption{Traffic network and aggregate measurements at initial and final times, respectively. The agents live on the edges and the thickness of the edges correspond to the population size of agents on the edges.}
	\label{fig:graph}
	\end{figure}

\begin{figure}[h]
	\centering
	\includegraphics[width=0.41\textwidth]{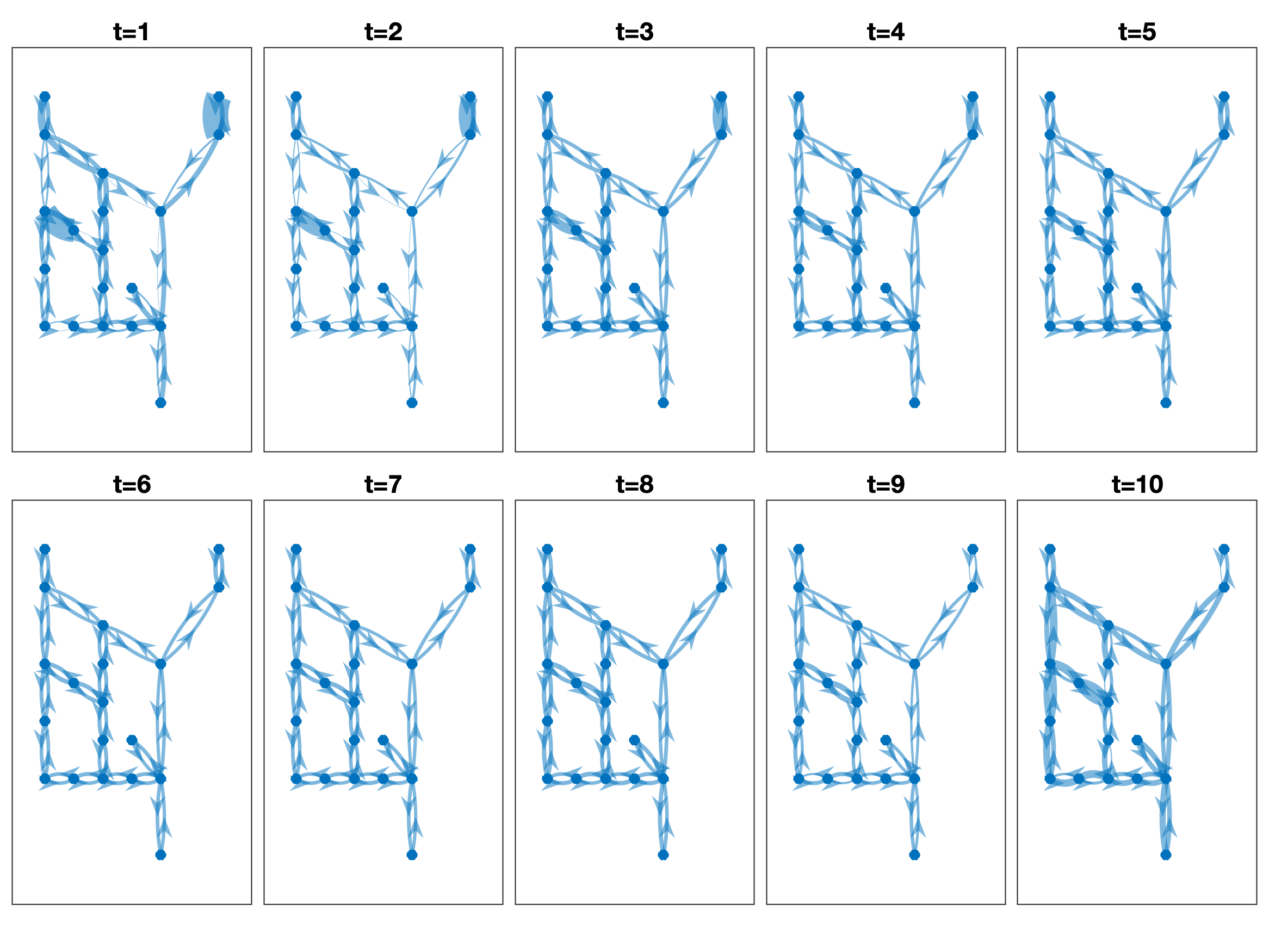}
	\caption{Estimated mass evolution of cars}
	\label{fig:car}
	\end{figure}

\begin{figure}[h]
	\centering
	\includegraphics[width=0.41\textwidth]{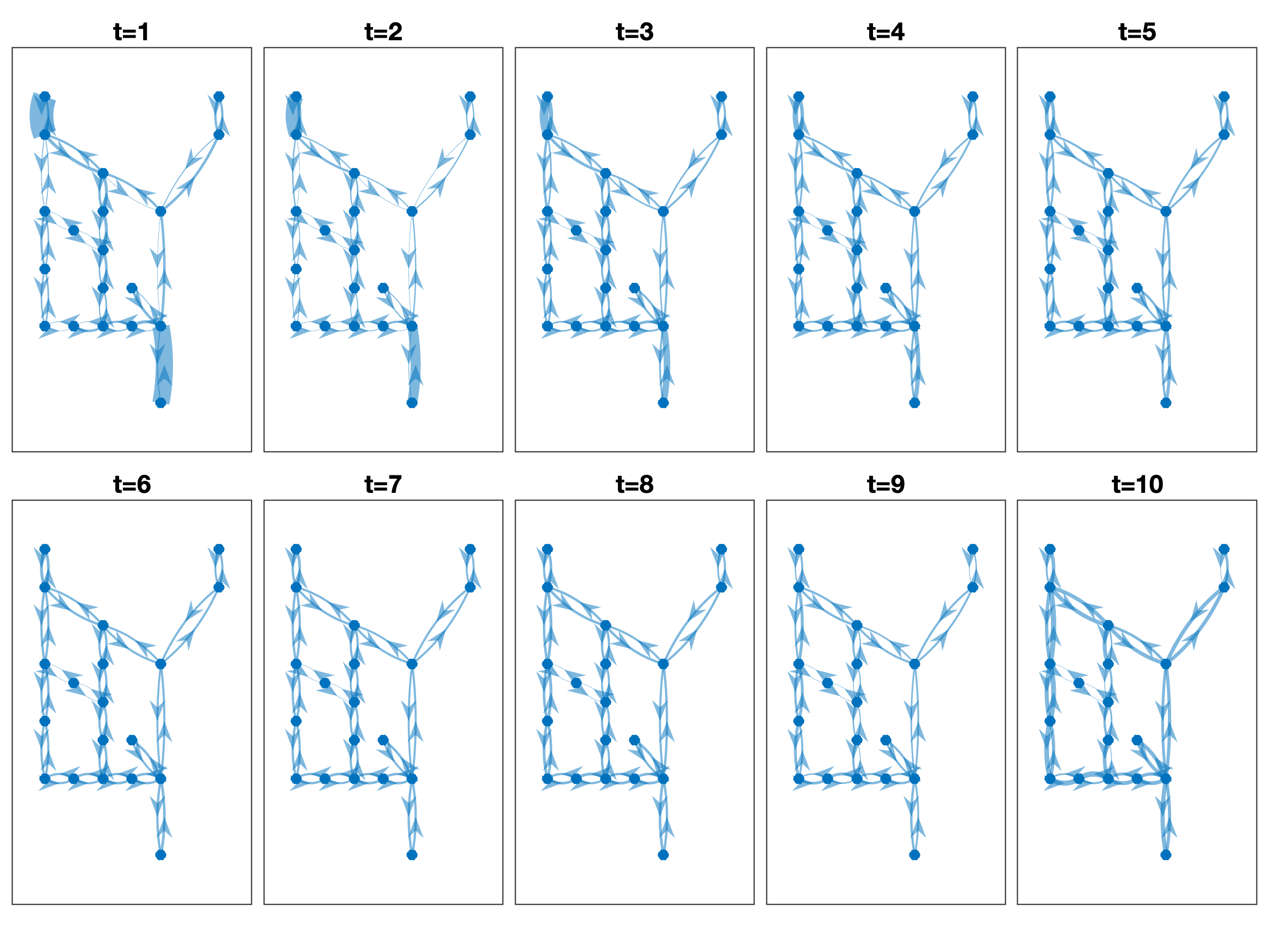}
	\caption{Estimated mass evolution of trucks}
	\label{fig:truck}
	\end{figure}

\begin{figure}[h]
	\centering
	\includegraphics[width=0.41\textwidth]{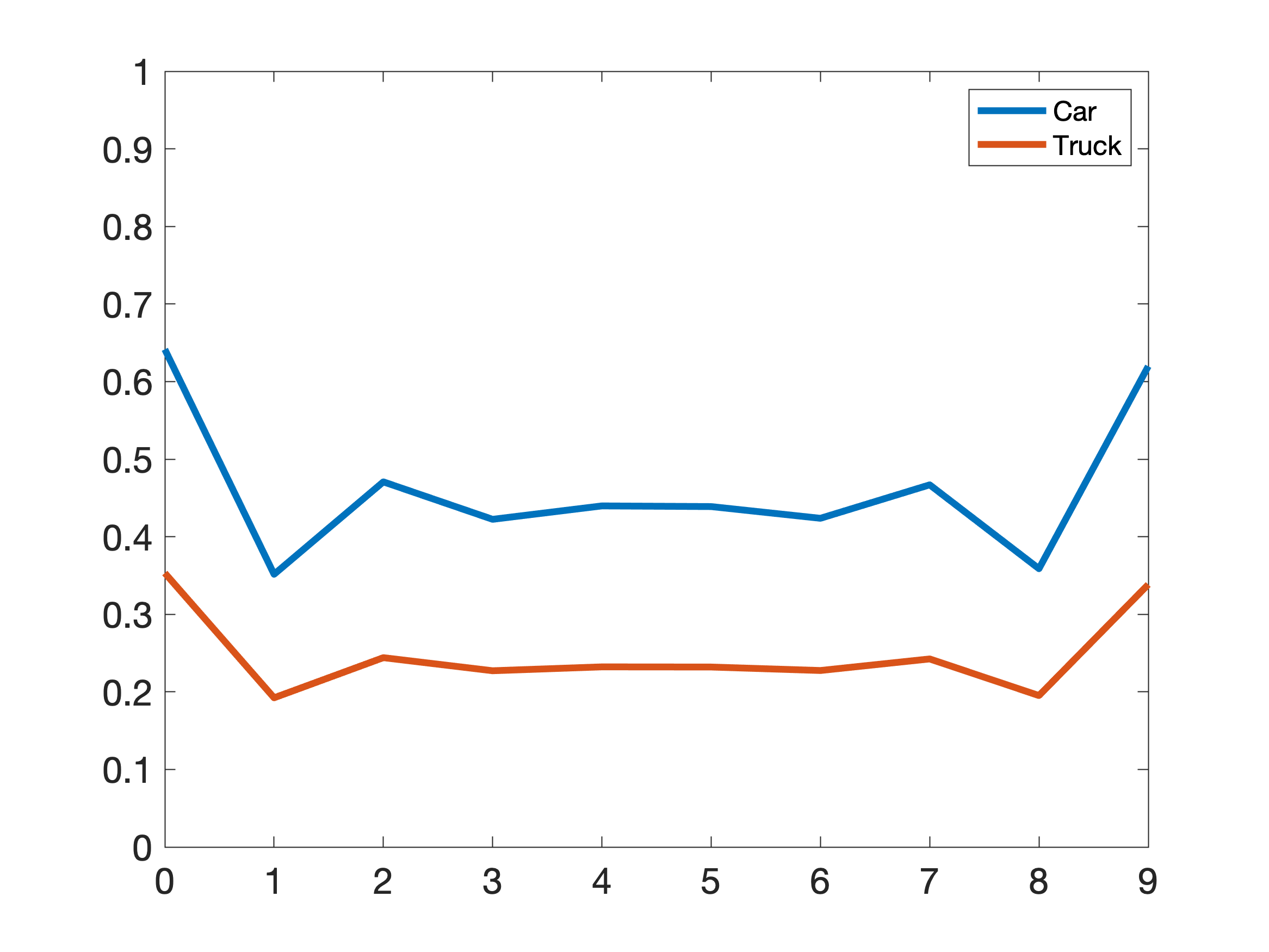}
	\caption{Evolutions of proportions of moving vechicles}
	\label{fig:mass}
	\end{figure}

In this section, we present a network flow example to illustrate our method. Consider a traffic network as shown in Figure \ref{fig:graph} where the agents live on the edges instead of nodes. Assume there are 2 commodities corresponding to cars and trucks. Assume the total number of individuals is $3000$: $2000$ cars and $1000$ trucks. Some nodes in the graph allow parking and thereby allow killing or creation with a certain probability. Each commodity follows certain prior dynamics that respect the topology of the network. In particular, the cars have a high probability to start from the edges $(2,\,4)$ and $(6,\, 9)$, and the trucks have a high probability to start from the edges $(1,\,3)$ and $(19,\, 18)$. Moreover, the trucks are not allowed on the edges $(6,\,9)$ and $(9,\, 10)$. We take as initial time $t=0$ and as a final time $t=9$.
In the traffic flow problem with aggregate observations, we assume the source and target distributions are given, as depicted in Figure \ref{fig:graph}. Our goal is to estimate the most likely evolution of each species,  or equivalently, to specify transition rates to effect the flow that reconciles the prior with the scarce data. We showcase the estimated results in Figures \ref{fig:car} and \ref{fig:truck}. These are deemed reasonable considering the limited aggregate measurements that are accessible 
and available, to reconcile with
their prior dynamics. We
display the result on
the proportion of cars and trucks that are not in the parking state in Figure \ref{fig:mass}. The transition, creation and killing rates of the posterior (solution to the problem) serve, both, as the solution to an estimation problem as well as a solution to a control problem, to dictate specifics of the flow (E.g., direction, parking or entering the flow by vehicles) that ensure matching the specified target marginals.

\section{Conclusions}

The purpose of the present work is to explain how the framework of Schr\"odinger bridges can be extended to model simultaneous transportation of multiple commodities, with partial aggregate data, as well as in the presence of creation and killing along the transport.

Traditionally, the term ``bridges'' refers to (probability laws on) paths linking marginal data, and the specific paradigm of Schr\"odinger bridges, is the method of constructing such laws on paths that maximize the likelihood over alternatives. To this end, a prior law is given (or chosen, as a ``design parameter'') and a posterior is sought to minimize the relative entropy between the two while maintaining consistency with available data. 
The first contribution of the paper is to note that, almost verbatim, the Schr\"odinger paradigm can be carried over to the case where multiple commodities are being transported at the same time and, when aggregate information is available on marginal distributions.

A second contribution is to enhance the Schr\"odinger framework with an additional state (or, states) that absorb commodities, and thereby take those out of the traffic flow randomly, or generate commodities similarly at random, so as to bring consistency with measured marginal distributions at various points. Here, marginals are assumed at the start and end of a specified interval.
Thus, this extra state (or, states) may be thought as a reservoir. 
The rate of killing in the prior is dictated by the distance to one of row sums in the transition kernel. The creation rate at nodes needs to be specified as a design parameter, or based on historical data.
Starting from such a suitably enlarged transition kernel that includes the reservoir state, the Schr\"odinger method can be readily applied to produce adjustment in the transition, creation, and killing probabilities so that the posterior kernel brings consistency with the measured marginal data in a way that may be deemed as the most likely.

The philosophy underlying the paper is fairly general, and can be suitably modified for more general marginal information that fully or partially reflect on the distributions of individual commodities. On a flip side, obtaining a transition kernel that meets marginals can be thought of as solving the control problem to decide on flow rates that effect an overall flow that matches the marginal data. In such a case marginal data can be seen as specifications dictated by supply and demand at various nodes in the network.

A final point that we wish to highlight is that the problems considered herein, and in the generality envisioned (though refrained so as to keep a simplified notation and exposition) can be efficiently cast and numerically solved as multimarginal transport problems.

{
\bibliographystyle{IEEEtran}
\bibliography{./refs}
}

\end{document}